\documentclass[10pt,a4paper]{amsart}

\usepackage{amssymb,amsfonts,latexsym}
\usepackage{xspace,url}
\usepackage{lscape}
\usepackage{mathtools}

\newtheorem{theorem}{Theorem}[section]
\newtheorem{prop}[theorem]{Proposition}

\newtheorem{cor}[theorem]{Corollary}
\newcommand{\Z}{\mathbb Z}

\long\def\Aut#1{{\rm Aut(#1)}}
\long\def\Inn#1{{\rm Inn(#1)}}
\long\def\Out#1{{\rm Out(#1)}}
\newcommand{\M}{{\rm M}}

\newcommand{\PGL}{{\rm PGL}}
\newcommand{\A}{{\rm A}}
\newcommand{\U}{{\rm U}}
\newcommand{\PSP}{{\rm PSP}}
\newcommand{\PgL}{{\rm P}\Gamma{\rm L}}

\newcommand{\SL}{{\rm SL}}
\newcommand{\PSL}{{\rm PSL}}
\newcommand{\ZC}{\textsf{(ZC-1)}\xspace}
\newcommand{\PQ}{\textsf{(PQ)}\xspace}
\newcommand{\SIPC}{\textsf{(SIP-C)}\xspace}
\newcommand{\pq}{$p \cdot q$\xspace}

\begin{document}

\begin{center} 
   \vbox to 2cm{}
   { \Huge \bf On the Gruenberg - Kegel Graph of Integral Group Rings of Finite Groups    
 }
\vskip 20pt
{ \large W.~Kimmerle\footnote{The 1st author has been supported by the DFG priority program SPP 1489}  and A.~Konovalov\footnote{The 2nd author has been supported by EP/M022641/1}} 
\end{center}

\begin{abstract}
 The prime graph question asks whether the Gruenberg - Kegel graph of
 an integral group ring $\Z G$ , i.e. the prime graph of the  
normalised unit group of $\Z G$ coincides with that one of the group $G$.
In this note we prove for finite groups $G$ a reduction of the prime
 graph question to almost simple groups. 
We apply this reduction to finite 
groups $G$ whose order is divisible by at most three primes and show
 that the Gruenberg - Kegel graph of such groups coincides with the
 prime graph of $G .$  
\end{abstract}

\vskip3em

\section{Introduction}

Let $G$ be a group. Its prime graph $\Pi (G)$ is defined as
follows. The vertices of $\Pi (G)$ are the primes $p$ for which $G$
has an element of order $p .$ Two different vertices $p$ and $q$ are joined
by an edge provided there is an element of $G$ of order \pq. If $R $
is a group ring then its Gruenberg - Kegel graph $\Gamma (R)$ is just the
prime graph $\Pi (V(R))$ of its normalised unit group $V(R) $.
\vskip1em
The integral group ring of $G$ is denoted by $\Z G$.
Its normalised unit group 
$V(\Z G) $ is the subgroup of the unit group $U(\Z G)$ which consists of 
all units with augmentation 1, called normalised units.
The objects of this article 
are $\Pi (V(\Z G))$ and $\Gamma (\Z G )$ 
in the case when $G$ is a finite group. 
\vskip1em
The question is whether $\Gamma (\Z G) $ coincides with
that of $G $ \cite[Problem 21]{JKNM}. 
This question, known as the ``Prime Graph Question'' (PQ) may be regarded as   
a weak version of the first Zassenhaus conjecture \ZC which says
that each torsion unit of $V(\Z G)$ is conjugate to an element of $G$,
where $G$ is considered in a natural way as subgroup of $V(\Z G)$ and
its elements are then called trivial units of $\Z G.$ The
conjecture \ZC is
certainly one of the major open questions for integral group rings and
if it is valid for a specific group $G$ it provides of course a positive answer
to the prime graph question for $V(\Z G) .$ However
only for specific classes (mainly specific soluble ones) of finite groups \ZC
is known. It is known that (PQ) holds for all soluble groups
\cite{KiPrimgraph}. In this 
article we show that (PQ) holds as well for large classes of finite insoluble
groups. 
\vskip1em
In Section 2 we reduce the study of the prime graph question to the study of 
nonabelian composition factors and their automorphism groups. 
Theorem \ref{reduce} shows that (PQ) holds provided $G$ does not map onto an
almost simple group which is a counterxample to (PQ). 
This gives 
many computer calculations made in the last years with respect to 
simple groups 
(see e.g. \cite{BJK,BKL,Hw4}) a new value. These calculations are not only examples. The verification of (PQ) 
for a single simple group establishes (PQ) for  infinitely many groups. 
Moreover with respect to sporadic simple groups the computational proof of (PQ) might be the only way. 
We note that not only simple groups have to be 
checked but also their automorphism groups have to be examined which has been 
done up to now only very rarely. 
\vskip1em
\vskip1em
As an application, we prove in Section 3 that for finite groups 
whose nonabelian 
composition factors have order divisible by only three primes the
prime graph question has a positive answer. 
\vskip1em
For this it is shown with
the aid of computational algebra that the
prime graph question, or even in some cases \ZC, has a positive
answer for all almost simple groups $G$ of this
type. For the most cases these computer calculations rest on two
independent implementations of the HeLP - algorithm (an algorithm initially
developed by I.~B.~S.~Passi and I.~Luthar, and later enhanced by M.~Hertweck).
One of them (ubpublished) is due to V.~Bovdi and the second author. The second 
one is due to A.~B\"achle and L.~Margolis and is available in the
{\sf HeLP} package \cite{help-package} which is redistributed together
with the computational algebra system GAP \cite{GAP}.
We like to point out that Theorem \ref{3primes} has been completed by a recent
result of A.~B\"achle and L.~Margolis for torsion units of integral group rings
of the groups $\M_{10}$  and $\PGL(2,9) \cong A6.2_2$ which is based
on a new algorithmic method, the so-called lattice method
\cite{Margthesis}. This method is also an essential tool in the
remarkable analysis of $\Gamma (\Z G)$ for almost simple groups $G$ whose order
is divisible by exactly four primes \cite{BaeMarFour1,BaeMarFour2}. The number
of such almost simple groups however is infinite and generic character
tables of series of almost simple groups enter the picture.
\vskip1em
With the aid of Theorem \ref{reduce} and generic character tables it follows that 
the prime graph question has a positive answer for all finite groups whose composition factors are either of
prime order or isomorphic to $PSL(2,p)$, where $p \geq 5$ is a prime (cf. Theorem \ref{psl}). 
    

\section{The reduction to almost simple groups}

An almost simple group of type $S$ is a subgroup of the automorphism group of a finite 
non-abelian simple group $S$ which contains $\Inn{S} \cong S .$ A
group $G$ is called almost simple if it is almost simple of type $S$
for some simple group $S .$ Clearly each finite simple group $S$
defines its own family of almost simple groups parametrized by the subgroups
of $\Out{S}$. It may happen that this
family just consists of $S .$ 
The object of this section is the following result which reduces the prime graph question to the study of almost simple groups.  
 
\vskip1em

\begin{theorem} \label{reduce}
Let $G$ be a finite group. Assume that for each almost simple group $X$ 
which occurs as image of $G$ the prime graph question has a positive answer.
Then the prime graph question has a positive answer for $G$.
\end{theorem}

\vskip1em

%
%
%
%

For $u \in V(\Z G)$ and $x \in G$ by $\varepsilon_x(u)$ we denote, as usual,
the partial augmentation of $u$ with respect to the conjugacy class $x^G$.

\begin{prop} \label{ext} 
Let $N$ be a normal subgroup of $E$ and $G = E/N$.
Assume that $\Pi(V(\Z G)) = \Pi (G)$.
Let $q \in \pi (E) \setminus \pi (N)$ and $p \in \pi (E)$ with $q \neq p.$ Then $V(\Z E)$ has elements of order
\pq if and only if $E$ has elements of this order. 
\end{prop}

{\bf Proof.} 
The proof in one direction is obvious. So 
assume that $V(\Z E)$ has a unit $u$ of order \pq but $E$ does
not have an element of this order.
If $\hat{\phi}(u)$ has order \pq,
by the assumption that $\Pi(V(\Z G)) = \Pi (G)$ we get 
an element of the order \pq in $G$. But then
$E$ has an element of this order.
It follows that $\hat{\phi} (u^p) = 1$ or $\hat{\phi}
(u^q) = 1$, where $\hat{\phi}$ denotes the homomorphism from $V(\Z E)$
to $V(\Z G)$ induced from $\phi .$

Let $v$ be a non-trivial torsion unit of $\Z E .$ 
By \cite[Theorem 2.7]{MRSW} we get that $\varepsilon_x( V) = 0 $ provided 
$o(x)$ 
is
divisible by a prime not dividing the order of $v .$ Moreover by
S.D.~Berman \cite{Ber} and
G.~Higman \cite{Hig,San} the 1-coefficient a of non-trivial torsion unit with augmentation 1 
has to be zero. It follows that orders of
torsion elements of prime order in
the kernel of $\hat{\phi} $ divide the order of $N.$ 
Because of the assumption that $q$ does not divide the order of $N$, we must have $\hat{\phi} (u^q) = 1 .$

For a group $X$ and a prime $r$ 
denote by $X(r)$ the set of all elements of $X$ which order is a positive 
power of $r$. 
Write $u$ as
$$ u = \sum_{g \in E(p)} a_gg + \sum_{g \in E(q)} b_gg + \sum_{g \in
  R} c_gg , $$
where $R = E \setminus (E(p) \cup E(q)) .$
Because $E$ has no elements of the order \pq, by \cite{MRSW} we conclude that
$\varepsilon_g(u) = 0 $ if $g \in R$.
From this it follows that $\sum_{g \in R} c_g = 0$.
Clearly

$$ \hat{\phi}(u) = 
\sum_{\mathclap{g \in E(p)}} a_g\phi(g) + 
\sum_{\mathclap{g \in E(q)}} b_g\phi(g) +  $$
$$ + \sum_{\mathclap{\substack{g \in R, \\ \phi(g) \in G(p) }}} c_g\phi(g) + 
\sum_{\mathclap{\substack{g \in R, \\ \phi(g) \in G(q)}}} c_g\phi(g) + 
\sum_{\mathclap{\substack{g \in R, \\ \phi(g) \in R'}}} c_g\phi(g), $$
where $R' = G \setminus (G(p) \cup G(q))$.
By the above, $\hat{\phi} (u)$ has order $q$, hence
$$\sum_{g \in E(p)} a_g \; + \;
\sum_{\mathclap{\substack{g \in R, \\ \phi(g) \in G(p) }}} c_g \; + \; 
\sum_{\mathclap{\substack{g \in R, \\ \phi(g) \in R'}}} c_g$$
is the sum of all partial augmentations of $\hat{\phi}(u)$
with respect to elements of $G$ whose order is not a positive power of $q$.
Because $\varepsilon_x(u)=0$ for each $x \in R$ and because $\phi $ maps conjugacy classes onto conjugacy
classes, we see that
$\sum_{g \in E(p)} a_g = 0$. 
Similarly we get $\sum_{g \in E(q)} b_g = 1$
using additionally the fact that the augmentation of
$\hat{\phi}(u)$ is equal to $1$.

Consider $$u^q = 
\sum_{g \in E(p)} \alpha_g g + 
\sum_{g \in E(q)} \beta_g g + 
\sum_{g \in R} \gamma_g g .$$
Because $u^q $ has order $p$, we see as before that $\sum_{g \in E(q)} \beta_g = 0$ and $\gamma_1=0$.

Denote as usual by $[\Z E,\Z E]$ the additive commutator of $\Z E$, i.e. the abelian subgroup generated by 
$\{ xy - yx \mid x, y \in E \}.$ Because $q $ is a prime, an easy calculation shows that

$$ u^q \equiv \sum_{g \in E( p)} a_g g^q  + \sum_{g \in E(q)} b_g g^q + \sum_{g \in R} c_g g^q \mod [\Z E,
\Z E] + q \Z E .$$ 

As $g^q$ is a $q$-element if and only if $g$ is a $q$-element, it follows that
$$\sum_{g \in E(q)} b_g \equiv \sum_{g \in E(q)} \beta_g \mod q.$$

Consequently,
$$ \sum_{g \in E(q)} b_g \equiv 0 \mod q .$$

This contradiction to ($*$) completes the proof. 
\qed
\vskip1em

\begin{cor} \label{complete} 
Let $G$ be a group and $N$ be a normal subgroup of $G$.
Then if $\Pi(N)$ is a complete graph and
$\Pi(V(\Z \; G/N)) = \Pi (G/N)$ then
$\Pi(V(\Z G)) = \Pi (G)$.
\end{cor}

{\bf Proof.} By Proposition \ref{ext}, we have only to consider
different primes dividing the order of $N$. 
By the assumption, these primes are connected in $\Pi (G)$ and therefore
also in $\Pi(V(\Z G))$.
\qed
\vskip1em

{\bf Remarks.} Instead of \cite[Theorem 2.7]{MRSW} we could have used
the stronger result \cite[Theorem 2.3]{Hw4} which says that partial
augmentations of a torsion unit of $\Z G$ are non-zero only for
conjugacy classes of $G$ whose representative has order dividing the
order of $u$. It is clear that this and modifications of the
arguments for Proposition \ref{ext} will lead to stronger results 
on the possible orders of torsion units in integral group
rings. However in this article we restrict ourselves on (PQ). \\

The special case when $N$ is a $p$-group in Corollary \ref{complete}
gives a proof of \cite[Proposition 4.3]{KiPrimgraph}. In the case of
soluble groups even more is known.
 
In \cite[Theorem]{Hw5} it is shown
that even \SIPC , i.e. that all finite cyclic subgroups of $V(\Z G)$ are
isomorphic to one of $G$, holds provided $G$ is
soluble. Note that \SIPC has been posted as research problem 8 in
\cite{Seh} and the abbreviation \SIPC stands for the subgroup isomorphism problem for cyclic subgroups.
A closer look shows that the proof given in \cite{Hw5} actually shows the
following.  

\begin{prop} \label{solubleext} Let $G$ be a finite group. Assume that $N$ is a soluble
  normal subgroup of $G$ and that for each torsion unit $u \in V(\Z G)$ there exists a $g \in G/N$ such that $o(u) = o(g) $ and $\varepsilon_g(u) \neq 0 $ then the same holds for torsion units of $\Z G .$ In particular  
  \SIPC is valid for $\Z G.$   
\end{prop}
\vskip1em

\begin{cor} Let $G$ be a finite Frobenius group. Then \SIPC is valid for $\Z G.$ More precisely 
for each torsion unit $u \in V(\Z G)$ there exists $g \in G$ such that $u$ and $g$ have the same order and
$\varepsilon_g(u) \neq 0 .$ 
\end{cor} 

\vskip1em
{\bf Proof.} From the structure theorems on Frobenius groups it follows that a finite Frobenbius group is  
soluble or a soluble extension of $S_5$ or $\SL(2,5) .$  In the first case the corollary follows from \cite[Theorem]{Hw5}, In the second case by \cite{DJM} for $\SL(2,5)$ and \cite{LuTra} for $S_5$  
the Zassenahus conjecture \ZC holds. This shows that the condition on the partial augmentations 
in Proposition \ref{solubleext} is fulfilled and the result follows. \qed 
\vskip1em
Proposition \ref{ext} permits also statements on some almost simple
groups.
\vskip1em

\begin{cor} Let $G$ be an almost simple group of type $S$. Assume that
  there is a prime $q \in G$ such that $q \notin \pi(S) .$ Then for
  each $p \in \pi (S)$ the primes $p$ and $q$ are connected in $\Gamma
  (\Z G)$ if and only if they are connected in $\Pi(G).$
\end{cor}

{\bf Proof.} By CFSG Schreier's  conjecture is valid, i.e. $\Out{S}$ is soluble. Thus $G/S$ is soluble and
therefore $\Gamma( \Z G/S) = \Pi (G/S).$ So we
may apply Proposition \ref{ext} and we get immediately the result. \qed 

\vskip1em

{\bf Proof of Theorem \ref{reduce}.} Clearly Theorem \ref{reduce} holds when $G$ itself is almost simple. 

Thus we may apply induction on the length of a chief series. 
If $G$ has a minimal normal subgroup $N$ which is an abelian $p$-group, then 
we see by Proposition \ref{ext} that the Theorem is valid provided it holds for $G/N .$ 

Suppose that all minimal normal subgroups of $G$ are perfect. Let 
$p,q \in \pi(G)$ be different primes and suppose that $V(\Z G)$ has an element
of order \pq. By Proposition \ref{ext} it follows that $G$ has an element of order \pq provided $p \in \pi (G) \setminus \pi (N)$ or 
$q \in \pi (G) \setminus \pi (N) $ for at least one minimal normal subgroup $N .$ So assume that $p$ and $q$ divide the order of each  minimal normal 
subgroup of $G .$ 

If $G$ has at least two different minimal normal subgroups $N_1, N_2$ 
then obviously $G$ has elements of order \pq because $N_1 \times N_2$ is a subgroup of $G .$ The same argument applies if $G$ has a minimal normal perfect subgroup which is not simple because such a subgroup is a direct 
product of at least two copies of isomorphic simple groups. 

Thus the only case which remains is that $G$ has a unique minimal normal 
subgroup which is a non-abelian simple group $S $ and $p$ and $q$ divide the order
of $S .$ But this case holds by assumption. \qed               

\vskip1em
Finally we give the following addendum to Theorem \ref{reduce}.
\vskip1em

\begin{prop} \label{2almost} Let $G$ be a finite group. Assume that $G$ has a normal subgroup $N$ such that $G/N$ is the direct product of two almost simple groups. Then (PQ) holds for $\Z G.$  
\end{prop}  

{\bf Proof.} For each normal subgroup $M \subset N$ the quotient $G/M$ satisfies the hypothesis of the proposition. 
Consider now a minimal counterexample and induct on a chief series of $G$ through $N .$  If $G/M$ maps onto $G/N$ then $G/M$ cannot be almost simple because $G/N$ is insoluble but by CFSG Schreier's conjecture holds and says that simple groups have soluble outer automorphism groups.  Thus a minimal counterexample is the direct product $G_1 \times G_2$ of two almost simple groups of type $S_1$, $S_2$ rsp. As in the proof of Theorem \ref{reduce} we see by Proposition \ref{ext} that primes $p$ and $q$ providing a counterexample to (PQ) both divide $|S_1|$ and $|S_2|.$ But then these primes are connected in $\Pi (S_1 \times S_2) $ and therefore in $\Pi (G).$ \qed      

\vskip1em

A typical example for Proposition \ref{2almost} is a wreath product of type $H \wr Q$, where $Q$ is the direct product of $A_5$ and $\PSL(2,7) .$  Theorem \ref{reduce} of course also shows that for finite groups which have no almost simple images, e.g. $A_5 \wr Q $ where $Q$ is an arbitrary soluble group, (PQ) has a positive answer.

\section{Groups of orders divisible by three primes only and further applications}\label{sect3primes}

The goal of this section is the following result.

\begin{theorem} \label{3primes} Let $G$ be a finite group such that the 
order of every almost simple image of $G$ is divisible by
exactly three different primes.
Then $\Pi(V(\Z G) = \Pi (G)$.
\end{theorem}

{\bf Proof.}
By Theorem \ref{reduce} it suffices to prove this for
all almost simple groups whose order is divisible by
exactly three different primes. 

By CFSG, 
the simple groups of order divisible by exactly three 
different primes are 
$$ PSL(2,5) \cong A_5, PSL(2,7), PSL(2,8), PSL(2,9) \cong A_6 $$
$$ PSL(2,17), PSL(3,3), PSP(3,4) \cong U(4,2), U(3,3) .$$
Their outer automorphism group has order 2, except the case
that $Out{A_6} \cong C_2 \times C_2$. 

In the table below we summarise details about these groups.
We also include in the table the isomorphism type of $\Out{S}$. 
In the third column we indicate with \ZC when the 1st Zassenhaus conjecture is established, with 
\SIPC when the finite cyclic subgroups of $V(\Z S)$ are isomorphic to subgroups of 
$S$ and with \PQ that the prime graph question has an affirmative answer. The fourth column contains the references.  

\vskip1em
\begin{center}

\begin{tabular}{|c|c|c|c|}\hline
                            & $\Out{S}$        &     & \\ \hline
$\A_5 \cong \PSL (2,5)$     & $C_2$            & \ZC & \cite{LuPas} \\
$\PSL(2,7)$                 & $C_2$            & \ZC & \cite{Hw3} \\
$\PSL(2,8)$                 & $C_3$            & \ZC & \cite{Gil,KiKo} \\
$\A_6 \cong \PSL (2,9)$     & $C_2 \times C_2$ & \ZC & \cite{HwA6}\\
$\PSL(2,17)$                & $C_2$            & \ZC & \cite{Gil} \\
$\PSL(3,3)$                 & $C_2$            & \PQ & \cite{BaeKi} \\
$\PSP (3,4) \cong \U (4,2)$ & $C_2$            & \PQ & \cite{KiKo} \\ 
$\U (3,3)$                  & $C_2$            & \SIPC & see below \\ \hline

\end{tabular}

\end{center}

\vskip1em

For the full automorphism groups of these simple groups the results are  as follows. 

\vskip1em

\begin{center}

\begin{tabular}{|c|c|c|} \hline
$\Aut{A_5} \cong S_5 $              & \ZC  & \cite{LuTra} \\
$\Aut{\PSL(2,7)} \cong \PGL (2,7)$   & \ZC  & \cite{KiKo} \\
$\Aut{\PSL(2,8)} \cong \PgL (2,8)$   & \SIPC & \cite{KiKo} \\
$\Aut{\A_6} \cong \PgL (2,9)$        & \SIPC & see below   \\
$\Aut{\PSL(2,17)} \cong \PGL (2,17)$ & \ZC  & \cite{KiKo} \\
$\Aut{\PSL(3,3)}$                    & \PQ  & \cite{KiKo} \\
$\Aut{\PSP (3,4)}$                   & \PQ  & \cite{KiKo} \\ 
$\Aut{\U (3,3)}$                     & \SIPC & see below   \\ \hline
\end{tabular}

\end{center}

For the three subgroups of index two of $\Aut{\A_6}$ the results are as follows. 

\vskip1em

\begin{center}

\begin{tabular}{|c|c|c|} \hline
$S_6$                               & \SIPC & \cite{KiKo} \\
$\PGL (2,9)$                        & \SIPC & \cite{BaeMar}\\
$M_{10}$                            & \SIPC & \cite{BaeMar}\\ \hline
\end{tabular}

\end{center}

The detailed report on calculations required to complete both
tables above was given in \cite{KiKo}; therefore, we will omit
them here and will only summarise known information in the table
below which presents results about possible orders and 
partial augmentations of normalised torsion units in integral group rings
for the groups listed above. 

\begin{itemize}
\item
Column 2 contains the name of the character table in the GAP character table
library provided by the GAP package {\sf CTblLib} \cite{ctbllib}.
\item 
Column 3 lists orders (of normalised torsion units) for which the rational
conjugacy is known, either as an immediate consequence of \cite[Proposition 3.1]{Hw3}
or using the HeLP--method.
\item
Column 4 lists orders of elements of $G$ with remaining non-trivial tuples 
of partial augmentations. In each entry of the form M(N), M means the order
and N means the number of all (both trivial and non-trivial) admissible 
tuples of partial augmentations that are produced by the
HeLP--method. 
\item
Column 5 lists orders of elements of $G$ which were omitted as not 
relevant to \PQ (for some groups it was possible to cover all or most of orders, though).
A dash (---) means that no orders were omitted.
\item
Column 6 lists orders that do not appear neither in $G$ nor in $V(\mathbb ZG)$,
eliminated using the HeLP--method (except units of order 6 in $A_6$ which are
eliminated in \cite{HwA6}, and of order 6 in $\PgL (2,9)$ and $M_{10}$ eliminated
in \cite{BaeMar}).
\item
Column 7 lists some orders with remaining non-trivial tuples of 
partial augmentations that still have to be eliminated for positive 
answers on \SIPC or \PQ for some groups.
\item
Column 8 lists further divisors of ${\rm exp}(G)$ needed
for the complete account on torsion units of $V(\mathbb ZG)$
but omitted as not relevant to \PQ.
\end{itemize}

\newpage
\begin{landscape}

\[
\footnotesize
\begin{array}{|c|c|c|c|c|c|c|c|}\hline
& & & & & & & \\
G & \text{Character table} & \text{\ZC} & \text{order(\#)}  & \text{Not considered} & \text{No orders}           &  \text{order(\#)} & \text{Not considered}  
\\
  & \text{name in {\sf CTblLib}}     &           & \text{in $G$}     & \text{orders in $G$}  & \text{in $V(\mathbb ZG)$}  & \text{in $V(\mathbb ZG)$} & \text{in $V(\mathbb ZG)$} 
\\ 
& & & & & & & \\
\hline
1 & 2 & 3 & 4 & 5 & 6 & 7 & 8 \\
\hline
PSL(2,5)
& \verb|A5|
& 2, 3, 5 
& \text{---} 
& \text{---} 
& 6, 10, 15 
& \text{---} 
&  \text{---} 
\\
PSL(2,7)
& \verb|PSL(2,7)| 
& 2, 3, 4, 7 
& \text{---} 
& \text{---} 
& 6, 14, 21 
& \text{---} 
&  \text{---} 
\\
PSL(2,8)
& \verb|PSL(2,8)| 
& 2, 3, 7, 9 
& \text{---} 
& \text{---} 
& 6, 14, 21 
& \text{---} 
&  \text{---} 
\\
A_6
& \verb|A6|
& 2, 3, 4, 5 
& \text{---} 
& \text{---} 
& 6, 10, 15 
& \text{---} 
&  \text{---} 
\\
PSL(2,17)
& \verb|PSL(2,17)| 
& 2, 3, 4, 8, 9, 17 
& \text{---} 
& \text{---} 
& 6, 34, 51 
& \text{---} 
&  \text{---} 
\\
PSL(3,3)
& \verb|PSL(3,3)| 
& 2, 4, 13 
& 3(5), 6(31) 
& 4, 6, 8 
& 26, 39 
& 12(4) 
&  \text{---} 
\\
U(4,2)
& \verb|U4(2)|
& 5 
& 2(3), 3(7), 4(13) 
& 6, 9, 12 
& 10, 15 
& \text{---} 
& 18 
\\
U(3,3)
& \verb|U3(3)|
& 2, 7 
& 3(3) 
& 4, 6, 8, 12 
& 14, 21, 24 
& \text{---} 
& \text{---} 
\\
\hline
\hline
S_5
& \verb|S5|
& 2,3,4,5,6 
& \text{---} 
& \text{---} 
& 10, 12, 15 
& \text{---} 
& \text{---} 
\\
PGL(2,7)
& \verb|L3(2).2|
& 2, 3, 4, 6, 7, 8 
& \text{---} 
& \text{---} 
& 12, 14, 21 
& \text{---} 
& \text{---} 
\\
\PgL (2,8)
& \verb|PSL(2,8).3| 
& 2, 3, 7, 9 
& 6(22) 
& \text{---} 
& 14, 18, 21 
& \text{---} 
& \text{---} 
\\
\PgL (2,9)
& \verb|A6.2^2|
& 2, 3, 5, 10 
& 4(5), 6(4), 8(3) 
& \text{---} 
& 12, 15, 20 
& \text{---} 
& \text{---} 
\\
\PgL (2,17)
& \verb|L2(17).2|
& 2, 3, 4, 6, 8, 9, 16, 17, 18 
& \text{---} 
& \text{---} 
& 12, 34, 51
& \text{---} 
& \text{---} 
\\
Aut{PSL(3,3)}
& \verb|L3(3).2|
& 2, 4, 8, 13 
& 3(5) 
& 6, 12 
& 26, 39 
& \text{---} 
& 24 
\\
\Aut{U(4,2)}
& \verb|U4(2).2|
& 5, 10 
& 2(10), 3(6), 9(2) 
& 4, 6, 8, 12
& 15 
& \text{---} 
& 18, 20, 24 
\\
\Aut{U(3,3)}
& \verb|U3(3).2|
& 2, 4, 7 
& 3(3), 6 (37), 8(3), 12(15) 
& \text{---}  
& 14, 21, 24 
& \text{---} 
& \text{---} 
\\
\hline
\hline
S_6
& \verb|S6|
& 3, 5 
& 2(4), 4(6), 6(16) 
& \text{---} 
& 10, 12, 15 
& \text{---} 
& \text{---} 
\\
\PgL (2,9)
& \verb|A6.2_2|
& 2, 3, 4, 5, 8, 10 
& \text{---} 
& \text{---} 
& 6, 15, 20 
&  \text{---} 
& \text{---} 
\\
M_{10}
& \verb|M10|
& 2, 3, 4, 5, 8 
& \text{---} 
& \text{---} 
& 6, 10, 15 
& \text{---} 
& \text{---} 
\\
\hline
\end{array}
\]
\end{landscape}

Let $G = \Aut{U(3,3)} = U(3,3).2$. To prove \SIPC for this group,
we need to show that there are no units of order 14, 21 and 24 in
$V(\Z G)$. Since $U(3,3) \subseteq U(3,3).2$, this will also prove
\SIPC for $U(3,3)$.

For torsion units of orders 14 and 21 we will use the method of
$(p,q)$-constant characters, introduced in \cite{BK}. Below we provide
the table containing the data for the constraints on partial augmentations,
from which the proof can be derived. The table uses the notation from \cite{BKL}.

\smallskip
\centerline{\small{
\begin{tabular}{|c|c|c|c|c|c|c|c|c|c|c|c|c|c|}
\hline
$|u|$&$p$&$q$& $\xi$                 & $\xi(C_p)$ & $\xi(C_q)$ & $l$ & $m_1$ & $m_p$ & $m_q$ \\
\hline
     &   &   &                       &            &            & 1   & 38    & 18    & 0     \\
14   & 2 & 7 & $\xi=(1,6,8)_{[3]}$   & 3          & 0          & 1   & 32    & 3     & 0     \\
     &   &   &                       &            &            & 7   & 32    & -18   & 0     \\
\hline
     &   &   &                       &            &            & 0   & 21    & 0     & -12   \\
21   & 3 & 7 & $\xi=(3,10)_{[*]}$    & 0          & -1         & 1   & 28    & 0     & -1    \\
     &   &   &                       &            &            & 7   & 21    & 0     & 6     \\
\hline
\end{tabular}
}}
\smallskip


%

For torsion units of order 24 we need to consider 39960 cases
determined by possible partial augmentations for torsion
units of order 2, 3, 4, 6, 8 and 12 (each having 2, 3, 3, 37, 4 and 15
cases respectively). Using computer calculations, in each of these cases we obtained
a contradiction by finding a pair of $i, j$ such that the
multiplicity of $\zeta^i$ as eigenvalue of $u$ in the 
representation affording $\chi_{j}$
is not an integer.

For the proof of Theorem \ref{3primes} it was a big advantage that there are 
only a finite number of finite simple groups whose order is divisible by
exactly three primes. If three is replaced by four then this is no longer
the case as mentioned in the introduction.
In the meantime there is also an infinite series of almost simple groups known
for which \PQ holds. Thus the following is a consequence from Theorem 
\ref{ext}.

\begin{cor} \label{psl}
Let $G$ be a finite group. Assume that all composition factors of $G$
are of prime order or isomorphic to $PSL(2,p)$ with a prime $p\geq 5$. 
Then PQ holds for $\Z G.$   
\end{cor}  

{\bf Proof.} By \cite[Prop 6.3, Prop 6.7]{Hw4} \PQ folds for all simple groups $PSL(2,p)$.
By \cite[Satz 2.4.2]{Margthesis} \PQ holds for all groups $PGL(2,p) =
\Aut{PSL(2,p)}$. Thus the result follows from Theorem \ref{ext}. \qed

\vskip1em

There are further simple groups $S$ known such that \PQ is valid for
$\Z G$ for all almost simple
groups $G$ of type $S ,$ e.g. when $S$ is isomorphic to one of the five
simple Mathieu groups or isomorphic to $A_7$ or $A_8 .$
This has been checked with HeLP package as well with the unpublished
implementation of the HeLP method by the 2nd author and V.~Bovdi. 

\vskip1em

\section{Reproducibility of results}

Since part of the results described in Section \ref{sect3primes}
are obtained with the help of computer calculations,
we are adding a special section to describe their
reproducibility. Following ACM policy \cite{acm-guidelines}, we distinguish
the ability to \emph{replicate} the experiment by 
obtaining same results with the same software on
a different computer from the ability to \emph{reproduce}
the experiment by obtaining same results using another
software.

In this terminology, we believe that our results are
reproducible. All calculations reported in the table
in Section \ref{sect3primes} were performed twice. One calculation
used the unpublished implementation of the Hertweck-Luthar-Passi
method developed by the 2nd author and Victor Bovdi and used, for
example, in \cite{BJK,BKL}. Another calculation used independent 
implementation of the HeLP method by A.~B\"achle 
and L.~Margolis provided in the HeLP package \cite{help-package},
which is now a part of the distribution for the computational
algebra system GAP. Both of these programs are written
in GAP, but developed independently, and use different
solvers: the former uses Minion \cite{Minion}, ECLiPSe \cite{eclipse}
and two custom solvers written in GAP, while the latter
uses Normaliz \cite{normaliz} by means of the GAP package 
NormalizInterface \cite{normaliz-gap}.

We compared both the final outcomes of calculations
and their numerical characteristics such as e.g. orders that have
to be analysed to reach certain conclusions, lists of 
admissible partial augmentations, etc., and detected no discrepancies;
hence we are convinced in the reproducibility of our results.

The check was tedious - sometimes it was not enough just to compare
lists of admissible partial augmentations straightforwardly. For
example, for PSL(3,3) it happened that HeLP package eliminates
the following cases of partial augmentations of elements of order 12:
$$ \{ ( -1, -1, -2, 1, 4 ), ( -1, 0, -3, 1, 4 ), ( -1, 0, 0, 1, 1 ), ( 1, 0, 0, 1, -1 ), $$
$$    ( 1, 0, 3, 1, -4 ), ( 1, 1, -1, 1, -1 ), ( 1, 1, 2, 1, -4 ) \},$$
which are not vanishing in calculations performed with the other implentation:
all of them have non-zero partial augmentations for the class of 
elements of order 2, while HeLP returns only solutions where the 
partial augmentations for this class are zero. It happened that
HeLP implements an additional constraint (so called ``Wagner test'',
see \cite{wagner-diplomarbeit}) unlike the other
code. Furthermore, as well as HeLP authors, we tried to use it to
check earlier results from e.g. \cite{BK}, and also they agree.

Of course, one could go further and observe that both programmes 
derive their data from the same source 
which is the GAP character table library \cite{ctbllib}. Here
\cite{rescience-atlas} discusses the reliability and reproducibility 
information contained in the Atlas of Finite Groups
and reports its verification done using Magma,
and \cite{rescience-tables} demonstrates how to compute in GAP the ordinary 
character tables of some Atlas groups, using
character theoretic methods.
This gives additional credibility and confidence in the 
correctness of our results.


\begin{thebibliography}{99}


\bibitem{acm-guidelines}
Association for Computing Machinery,
Result and Artifact Review and Badging
(\url{http://www.acm.org/publications/policies/artifact-review-badging}).

\bibitem{BaeKi} 
A.~B\"achle and W.~Kimmerle, 
On torsion subgroups in integral group rings of finite groups, 
\emph{J. Algebra} {\bf 326} (2011), 34--46.  

\bibitem{BaeMar}
A.~B\"achle and L.~Margolis, 
Rational conjugacy of torsion units in integral group rings of non-solvable groups,
to appear in Proc. Edinburgh Math.Soc.2017, arXiv:1305.7419v2 [math.RT], February 2014. 

\bibitem{help-package}
A.~B\"achle and L.~Margolis,
\emph{HeLP -- Hertweck-Luthar-Passi method, Version 3.0}; GAP package, 2016
(\url{http://homepages.vub.ac.be/abachle/help/}).

\bibitem{BaeMarFour1} {\rm A.~B\"achle and L.~Margolis},
On the Prime Graph Question for Integral Group Rings of 4-primary groups I,
ArXiv:1601.05689 [math.RT], January 2016. 

\bibitem{BaeMarFour2} {\rm A.~B\"achle and L.~Margolis},
On the Prime Graph Question for Integral Group Rings of 4-primary groups II,
ArXiv:1606.01506 [math.RT], June 2016. 

\bibitem{Ber}
{\rm S.~D.~Berman}, On certain properties of integral group rings (Russian), 
Dokl. Akad. Nauk SSSR (N.S.) {\bf 91} (1953), 7 -- 9.  


\bibitem{BK}
{\rm V.~Bovdi, A.~Konovalov}, 
Torsion units in integral group rings of Higman-Sims simple group, 
Studia Sci. Math. Hungar. {\bf 47} (2010), no.1, 1 -- 11.  

\bibitem{BJK}
{\rm V.~Bovdi, E.~Jespers, A.~Konovalov}, 
Torsion units in integral group rings of Janko simple groups, 
Math.Comp. {\bf 80} (2011), no.273, 593 -- 615.  

\bibitem{BKL}
{\rm V.~Bovdi, A.~Konovalov, S.~Linton},
Torsion units in integral group rings of Conway simple groups,
Int. J. Algebra Comput. 21 (2011), no.4, 615 -- 634. 


\bibitem{ctbllib}
{\rm T.~Breuer},
\emph{The GAP Character Table Library, Version 1.2.2}; GAP package, 2013,
\url{http://www.math.rwth-aachen.de/~Thomas.Breuer/ctbllib}.

\bibitem{rescience-atlas} 
{\rm T.~Breuer, G.~Malle, E.~A.~O'Brien},
Reliability and reproducibility of Atlas information,
arXiv:1603.08650 [math.GR], March 2016. 

\bibitem{rescience-tables} 
{\rm T.~Breuer},
Constructing the ordinary character tables of some Atlas groups using character theoretic methods,
arXiv:1604.00754 [math.RT], April 2016. 

\bibitem{normaliz}
{\rm W.~Bruns, B.~Ichim, T.~R\"omer, R.~Sieg and C.~S\"oger}: 
Normaliz. Algorithms for rational cones and affine monoids. 
Available at \url{https://www.normaliz.uni-osnabrueck.de}.


\bibitem{DJM}
{\rm M.~Dokuchaev, S.~O.~Juriaans and C.Polcino Milies}, Integral Group Rings of Frobenius Groups and the Conjectures of
H.J.Zassenhaus, Comm. in Alg. {\bf 25}(7), 2311 -- 2325 (1997).

\bibitem{eclipse}
{\rm J.~Schimpf and K.~Shen},
ECLiPSe - from LP to CLP,
Theory and Practice of Logic Programming, vol.12,
Special Issue on Prolog Systems 1-2, 127--156. 

\bibitem{GAP} 
{\rm The~GAP~Group}, \emph{GAP --- Groups, Algorithms, and Programming, 
Version 4.8.4}; 2016, \url{http://www.gap-system.org}.

\bibitem{Minion}
{\rm I.P.~Gent, C.~Jefferson and I.~Miguel},
MINION: A Fast, Scalable, Constraint Solver,
in Proceedings of the 17th European Conference on Artificial Intelligence (ECAI 2006).

\bibitem{Gil}
J.~Gildea, 
Zassenhaus conjecture for integral group ring of simple linear groups, 
\emph{J. Alg. Appl.}, \textbf{12} (2013), 1350016.

\bibitem{normaliz-gap}
{\rm S.~Gutsche, M.~Horn, C.~S\"oger}: 
NormalizInterface, GAP wrapper for Normaliz,
Version 0.9.8 (2016), (GAP package),
\url{https://gap-packages.github.io/NormalizInterface}.






\bibitem{Hw3}
{\rm M.~Hertweck}, On the torsion units of some integral group rings, 
Algebra Colloq. {\bf 13}(2) (2006) 329--348.

\bibitem{Hw4} {\rm M.~Hertweck}, 
Partial augmentations and Brauer character values of torsion
  units in group rings, e-print arXiv:math.RA/0612429v2.

\bibitem{Hw5}
{\rm M.~Hertweck},
The orders of torsion units in integral group rings of fiite solvable
groups, Comm.Alg. {\bf 36}, (2008) 3585 -- 3588.


\bibitem{HwA6}
M.~Hertweck, 
Zassenhaus Conjecture for $A_6$, 
\emph{Proc. Indian Acad. Sci. (Math.Sci)} 
\textbf{118} (2008), no.2, 189--195.  




\bibitem{Hig}
{\rm G.~Higman}, Units in group rings, D.Phil.thesis, Oxford Univ. 1940. 






\bibitem{JKNM} {\rm E.~Jespers, W.~Kimmerle, G.~Nebe and Z.~Marciciak}, Miniworkshop 
Arithmetik von Gruppenringen, vol.~4 Oberwolfach Reports
No.55, pp. 3149 -- 3179, EMS, Z\"urich 2007.


\bibitem{KiPrimgraph} {\rm W.~Kimmerle},  
On the prime graph of the unit group of integral group rings of finite groups, Contemporary Mathematics Vol.{\bf 420} 2006, 215 -- 228. 
  

\bibitem{KiKo}
W.~Kimmerle and A.~Konovalov,
On the prime graph of the unit group of integral group rings of finite
groups II, \emph{Stuttgarter Mathematische Berichte} 2012--018 
(\url{http://www.mathematik.uni-stuttgart.de/fachbereich/math-berichte/listen.jsp}). 





\bibitem{LuPas}
{\rm I.~S.~Luthar and I.~B.~S.~Passi}, Zassenhaus conjecture for $A_5$, 
J. Nat. Acad. Math. India {\bf 99} (1989), 1 -- 5.


\bibitem{LuTra}
{\rm I.~S.~Luthar and P.~Trama}, Zassenhaus conjecture for $S_5$, Comm. 
Alg. {\bf 19(8)} (1991), 2353 -- 2362 

\bibitem{Margthesis}
{\rm L.~Margolis}, 
Torsionseinheiten in ganzzahligen Gruppenringen nicht aufl\"osbarer Gruppen.
PhD thesis, Stuttgart, 2015. 

\bibitem{MRSW}
{\rm Z.~S.~Marciniak, J.~Ritter, S.~K.~Sehgal and A.~Weiss}, Torsion Units
in Integral Group Rings of Some Metabelian Groups, II, J. Number Theory 
{\bf 25} (1987), 340 -- 352.  




  

 
\bibitem{San} 
{\rm R.~Sandling},
Graham Higman's thesis ``Units in Group Rings" in
Integral Representations and Applications ed.  K.W.~Roggenkamp,
Springer Lecture Notes {\bf 882}  93--116, 1981. 


\bibitem{Seh}
{\rm S.~K.~Sehgal}, Units in integral group rings, 
Pitman Monographs and Surveys in
Pure and Applied Mathematics, Longman Scientific \& Technical, Essex, 1993.



 



\bibitem{wagner-diplomarbeit}
{\rm R.~Wagner},
Zassenhausvermutung \"uber die Gruppen $\PSL(2, p)$,
Diplomarbeit, Universit\"at Stuttgart, 1995.




\end{thebibliography}
\end{document}